\documentclass[11pt]{article}
\usepackage{geometry,amssymb,latexsym,amsmath,graphicx, enumerate}
\usepackage[utf8]{inputenc}

\usepackage{graphicx}
\usepackage{tikz,ifthen}
\usetikzlibrary{calc}

\newcommand{\qed}{\hfill ~$\square$\bigskip}
\newcommand{\proof}{\noindent{\bf Proof.} }

\newcommand{\D}{{Dominator}}
\newcommand{\St}{{Staller}}

\newtheorem{theorem}{Theorem}%[section]
\newtheorem{lemma}[theorem]{Lemma}
\newtheorem{corollary}[theorem]{Corollary}

\newtheorem{observation}[theorem]{Observation}

\newtheorem{question}{Question}
\begin{document}

\title{Complexity of the Game Domination Problem}

\author{
Bo\v{s}tjan Bre\v{s}ar $^{a,b}$
\and
Paul Dorbec $^{c,d}$
\and
Sandi Klav\v zar $^{e,a,b}$
\and
Ga\v{s}per Ko\v{s}mrlj $^{e}$
\and
Gabriel Renault $^{c,d}$
}

\date{}

\maketitle

\begin{center}
$^a$ Faculty of Natural Sciences and Mathematics, University of Maribor, Slovenia\\
\medskip

$^b$ Institute of Mathematics, Physics and Mechanics, Ljubljana, Slovenia\\
\medskip

$^c$ Univ. Bordeaux, LaBRI, UMR 5800, F-33400 Talence, France\\
\medskip

$^d$ CNRS, LaBRI, UMR 5800, F-33400 Talence, France\\
\medskip

$^e$ Faculty of Mathematics and Physics, University of Ljubljana, Slovenia\\
\end{center}

\begin{abstract}
The game domination number is a graph invariant that arises from a game,
which is related to graph domination in a similar way as the
game chromatic number is related to graph coloring. In this paper we show
that deciding whether the game domination number of a graph is bounded
by a given integer is PSPACE-complete. This contrasts the situation of
the game coloring problem whose complexity is still unknown.
\end{abstract}


\medskip\noindent
{\bf AMS Subj. Class:} 05C57, 68Q15, 05C69
\footnotetext[1]{\begin{minipage}[t]{.55\textwidth}© 2016. This manuscript version is made available under the CC-BY-NC-ND 4.0 license
\end{minipage} \hfill \raisebox{-4ex}{\includegraphics[height=1cm]{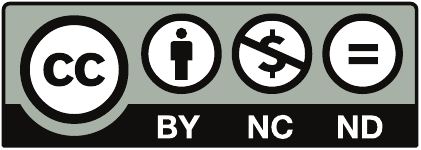}}}

%%%%%%%%%%%%%%%%%%%%%%%%%%%%%%%%%%%%%%%%%%%%%%%%%%%%%%%%%%%%%%%%%
\section{Introduction}
%%%%%%%%%%%%%%%%%%%%%%%%%%%%%%%%%%%%%%%%%%%%%%%%%%%%%%%%%%%%%%%%%

In a graph $G$, we denote by $N[v]$ the closed neighborhood of $v$, i.e. the vertex $v$ itself together with all its neighbors. By extension, the closed neighborhood $N[D]$ of a subset of vertices is defined by $N[D]=\cup_{v\in D} N[v]$.
%%%%% We need it for the formal definition later on  %%%%
 We say that a vertex $u$ (resp.\ a set of vertices $D$) {\em dominates} a vertex $x$ (and that $x$ {\em is dominated} by $u$ (resp.\ $D$)), if $x\in N[u]$ (resp.\ $x\in N[D]$). If a set of vertices $D$ dominates all the vertices of $G$, we say that $D$ is a {\em dominating set} of $G$. The minimum cardinality of a dominating set of a graph $G$ is called the {\em domination number} of $G$, and is denoted by $\gamma(G)$.   

The domination game, introduced in~\cite{brklra-2010}, is played by
two players on an arbitrary graph $G$. The two players are called \D\
and \St, indicating the roles they are supposed to play in the
game. They are taking turns choosing a vertex from $G$ such that
whenever they choose a vertex, it dominates some vertex not dominated by the vertices chosen in earlier turns.
Assuming that \D\ starts the game, and letting $d_1, d_2, \ldots$ denote the sequence of vertices chosen by \D\ and $s_1,s_2, \ldots$ the sequence chosen by Staller, the condition the two players must fulfill is that for each $i$,
\begin{itemize}
\item $N[d_i] \setminus \cup_{j=1}^{i-1}N[\{d_j, s_j\}]\not=\emptyset$; and
\item $N[s_i] \setminus \left(\cup_{j=1}^{i-1}N[\{d_j, s_j\}]  \cup N[d_i] \right)\not= \emptyset$.
\end{itemize}

\noindent The game ends when all vertices of $G$ are dominated.
The aim of \D\ is that the total number of moves played in the game
is as small as possible, while \St\ wishes to maximize this number.
By {\em D-Game} we mean a game in which Dominator has the first
move, while {\em S-Game} refers to a game in which \St\ begins.
Assuming that both players play optimally, the {\em game domination
number} $\gamma_g(G)$ of a graph $G$ denotes the number of moves played in D-Game, and the {\em Staller-start game domination number} $\gamma_g'(G)$ the number of moves played in S-Game

Giving a formula for the exact value of $\gamma_g(G)$ and $\gamma'_g(G)$ is usually a difficult problem, and is resolved only for some very simple families of graphs $G$, e.g.\ for paths and cycles~\cite{bill-2012+} and for combs~\cite{gk-14}. It was conjectured in~\cite{bill-2013} that the upper bound $\gamma_g(G)\le 3/5|V(G)|$ holds for any isolate-free forest as well as for any isolate-free graph $G$. In the seminal paper~\cite{bill-2013}, the conjecture was verified for forests in which each component is a caterpillar. Bujt\'{a}s followed with a powerful discharging-like method and proved the conjecture for forests in which no two leaves are at distance~$4$ apart~\cite{bu-2013, bujtas-2015a}. Using the method, Schmidt~\cite{sc-2016+} determined a larger class of forests for which the conjecture holds. Bujt\'{a}s~\cite{bujtas-2015b} proved the conjecture for graphs with the minimum degree at least $3$, while Henning and Kinnersley~\cite{heki-2016} further extended the truth of the conjecture to graphs of minimum degree at least $2$. Note also that in~\cite{brklko-2013}, large families of trees were constructed that attain the conjectured $3/5$-bound and all extremal trees on up to 20 vertices were found.

The game domination number has been studied from several additional aspects, see~\cite{brdokl-2014, brklra-2013, buklko-2015, doko-2015, klkosc-2016, nasi-2016+}. On the other hand, the only contribution so far from the  algorithmic point of view is the paper~\cite{klkosc-2015}. This comes with no surprise, considering that the complexity of the much older coloring game~\cite{bagrkizh-07} has not been determined yet. Bodlaender~\cite{boa-91} proved that a version of the coloring game where the order of the vertices to be colored is prescribed in advance is PSPACE-complete  (see also~\cite{bobuchjekr-09}), but to the best of our knowledge, no result is known for the standard coloring game.

In~\cite{klkosc-2015} it is shown that for a given integer $m$ and a given graph $G$, deciding whether $\gamma_g(G)\le m$ can be done in $\mathcal O(\Delta(G)\cdot |V(G)|^m)$ time. 
% This means that if $k$ is {\em not} part of the input, the problem becomes polynomial.
In this paper we complement this result by proving that the complexity of verifying whether the game domination number of a graph is bounded by a given integer is in the class of PSPACE-complete problems, implying that every problem solvable in polynomial space (possibly with exponential time) can be reduced to this problem. In particular, this shows that the game domination number of a graph is harder to compute than any other classical domination parameter (which are generally NP-hard), unless NP=PSPACE. The reduction we use can be computed with a working space of logarithmic size with respect to the entry, making this problem $\log$-complete in PSPACE. (For additional problems that were recently proved to be PSPACE-complete see~\cite{atserian-2014,fomin-2014,grier-2013,ito-2014}.)

 %In particular, if $S=\{x\}$ we write $G|x$.
%With these concepts in hand we now recall the following fundamental tool due to Kinnersley, West, and Zamani~\cite{bill-2013}.

%\begin{lemma}
%{\rm (Continuation Principle)}
%{\rm\cite[Lemma 2.1]{bill-2013}}
%\label{lem:continuation}
%Let $G$ be a graph and $A, B\subseteq V(G)$.
%If $B\subseteq A$, then $\gamma_g(G|A) \le \gamma_g(G|B)$ and
%$\gamma_g'(G|A) \le \gamma_g'(G|B)$.
%\end{lemma}

%We also consider a game (with respect to some subgraph of the whole graph in which the game is played) in which one of the players passes a move.

In the following section, we present a reduction from the classical PSPACE-complete problem POS-CNF to a game domination problem where some vertices are set to be dominated before the game begins. Then, in Section~\ref{sec:conclude}, we describe how to extend the reduction to the Staller-start domination game and to the game on a graph not partially-dominated. We conclude with some open questions.

%%%%%%%%%%%%%%%%%%%%%%%%%%%%%%%%%%%%%%%%%%%%%%%%%%%%%%%%%%%%%%%%%
\section{PSPACE complexity of the game domination problem}
\label{sec:main}
%%%%%%%%%%%%%%%%%%%%%%%%%%%%%%%%%%%%%%%%%%%%%%%%%%%%%%%%%%%%%%%%%

The game domination problem is the following:

\begin{center}
\fbox{\parbox{0.85\linewidth}{\noindent
{\sc Game Domination Problem}\\[.8ex]
\begin{tabular*}{.93\textwidth}{rl}
{\em Input:} & A graph $G$, and an integer $m$.\\
{\em Question:} & Is $\gamma_{g}(G)\le m$?
\end{tabular*}
}}
\end{center}

To prove the complexity of the {\sc Game Domination Problem}, we propose a reduction from the POS-CNF problem, which is known to be log-complete in PSPACE~\cite{ts-78}.
In POS-CNF we are given a set of $k$ variables, and a formula that is a conjunction of $n$ disjunctive clauses, in which only positive variables appear (that is, no negations of variables). Two players alternate turns, Player 1 setting a previously unset variable TRUE, and Player~2 setting one FALSE. After all $k$ variables are set, Player 1 wins if the formula is TRUE, otherwise Player~2 wins. 

In the proof of our main result, we will make use of the so-called {\em partially-dominated graph}. This is a graph together with a declaration that some vertices are already dominated, that is, they need not be dominated in the rest of the game. For a vertex subset $S$ of a graph $G$, let $G|S$ denote the partially dominated graph
in which vertices from $S$ are already dominated (note that $S$ can be an arbitrary subset of $V(G)$, and not only a union of closed neighborhoods of some vertices).
We transform a given formula ${\cal F}$ using $k$ variables and $n$ disjunctive clauses into a partially dominated graph $G_{\cal F}|A$, having $9k+n+4$ vertices.
We then prove that Player 1 has a winning strategy for a formula ${\cal F}$ if and only if $\gamma_g(G_{\cal F|A})\le 3k+2$.

\begin{figure}[htb!]
\begin{center}

\begin{tikzpicture}
[thick,scale=1,
     vertex/.style={circle,draw,inner sep=0pt,minimum
size=1.5mm,fill=white!100},
     blackvertex/.style={circle,draw,inner sep=0pt,minimum
size=1.5mm,fill=black!100},
     clause/.style={circle,draw,inner sep=0pt,minimum
size=3mm,fill=white!100}]

% positions

%\coordinate (c) at (0,-2);
\coordinate (d1) at (0,0);
\coordinate (d2) at (0,3.5);
\coordinate (s1) at (-2,1);
\coordinate (s2) at (2,1);
\coordinate (x0) at (0,1);
\coordinate (x11) at (-1,1.5);
\coordinate (x12) at (-1,2.5);
\coordinate (x21) at (1,1.5);
\coordinate (x22) at (1,2.5);

%edges

\draw (d1)--(s1)--(x12)--(d2)--(x22)--(s2)--(d1);
\draw (d1)--(x0);
\draw (s2)--(x21)--(x0)--(x11)--(s1);
\draw
(x12)--(x0)--(x11)--(x22)--(x12)--(x21)--(d2)--(x11)--(x21)--(x0)--(x22);

%\draw (c)--(d1);
%\draw (c) arc (-90:90:2.775);

%vertices
%\draw (c) node[clause] {};
\draw (d1) node[blackvertex] {};
\draw (d2) node[blackvertex] {};
\draw (s1) node[vertex] {};
\draw (s2) node[vertex] {};
\draw (x0) node[vertex] {};
\draw (x11) node[vertex] {};
\draw (x12) node[vertex] {};
\draw (x21) node[vertex] {};
\draw (x22) node[vertex] {};

%labels
%\draw (c) node[left] at +(-0.2,0) {$c$};
\draw (d1) node[below right] {$a_1$};
\draw (d2) node[above right] {$a_2$};
\draw (s1) node[left] {$b_1$};
\draw (s2) node[right] {$b_2$};
\draw (x0) node[below right] {$z$};
\draw (x11) node[below] {$x$};
\draw (x12) node[left] {$x'$};
\draw (x21) node[below] {$y$};
\draw (x22) node[right] {$y'$};

\end{tikzpicture}

\end{center}
\caption{The widget graph $W$ used to represent each variable}
\label{fig:variable}
\end{figure}
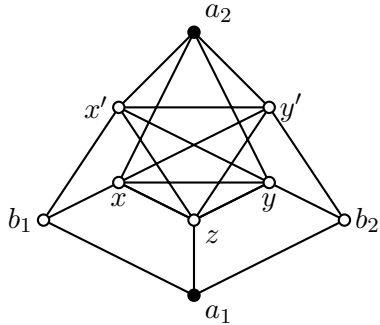

In the construction of $G_{\cal F}$, we use $k$ copies of the widget graph $W$ in correspondence with the $k$ variables. The graph $W$ is constructed from the disjoint union of the cocktail-party graph $K_6-M$ on the vertex set $\{a_2,x,x',y,y',z\}$ with $a_2z,xx',yy'\not\in E(W)$,
and of the path $P:b_1a_1b_2$, by the addition of the edges $b_1x,b_1x',b_2y,b_2y'$ and $a_1z$.
Moreover, the vertices $a_1$ and $a_2$ are assumed to be dominated, that is, we are considering the partially dominated graph $W|\{a_1,a_2\}$. In Fig.~\ref{fig:variable} the graph $W$ is shown, where the vertices $a_1$ and $a_2$ are filled black to indicate that they are assumed to be dominated.

Although the rules of the game do not allow any of the players passing a move, it may happen during the game played in $G_{\cal F}$ that after some moves were played in a subgraph $W_X$ isomorphic to the widget graph $W$ but before all vertices of this subgraph are dominated, one of the players chooses to make a move outside $W_X$. In this case, the other player could be the first one to play again in the subgraph $W_X$. From the point of view of the graph $W$, this corresponds to one player passing a move in $W$.

For a better distinction between the two games, in the S-game we denote $s_1',s_2', \ldots$ the sequence of vertices chosen by \St\ and $d_1',d_2', \ldots $ those chosen by \D, while in the D-game we use notation $d_1, d_2, \ldots$ for the vertices chosen by \D\ and $s_1,s_2, \ldots$ for \St's moves.
In the proof of the main result, we use several times the following properties of the graph $W$.

\begin{observation}
\label{obs:widget}
If $W$ is the
 widget graph, then
\begin{enumerate}[(i)]
\item
$\gamma_g(W|\{a_1,a_2\})=3$ and $a_1$ is an optimal first move;
\item
if $d_1=a_1$, and \St\ passes her first move, then \D\ can finish the game in $W|\{a_1,a_2\}$ in two moves (by playing $a_2$);
\item
if $d_1=a_1$ (or $d_1=a_2$), $s_1=b_1$, and \D\ passes his second move, then \St\ can ensure four moves will be played in $W|\{a_1,a_2\}$;
\item
$\gamma_g'(W|\{a_1,a_2\})=3$ and $b_1$ is an optimal first move for \St\ (an optimal response to $s_1'=b_1$ is $d_1'=z$);
\item
if $s'_1=b_1$ and \D\ responds playing $a_1$ or $a_2$, then \St\ can enforce four moves are played in $W|\{a_1,a_2\}$ (by playing respectively $y$ or $x$);
\item
if \St\ starts and \D\ passes his first move, then after any second move of \St, \D\ can finish the game in $W|\{a_1,a_2\}$ with the third move, ensuring in addition that $a_1$ or $a_2$ is played.
\qed
\end{enumerate}
\end{observation}

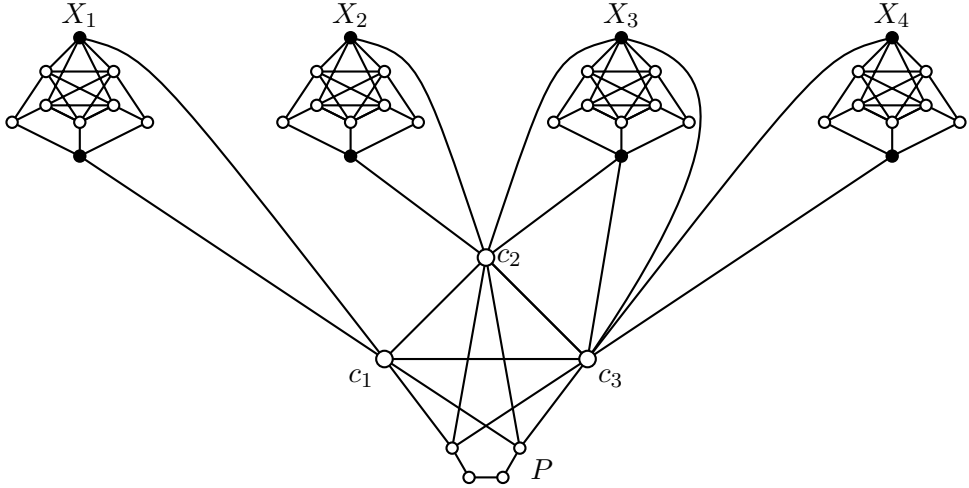
\begin{figure}[htb!]
\begin{center}

\begin{tikzpicture}
[thick,scale=0.45,
     vertex/.style={circle,draw,inner sep=0pt,minimum
size=1.5mm,fill=white!100},
     blackvertex/.style={circle,draw,inner sep=0pt,minimum
size=1.5mm,fill=black!100},
     clause/.style={circle,draw,inner sep=0pt,minimum
size=2.2mm,fill=white!100}]

\def\k {3};
\def\l {4};

%%%% Variables

\foreach \i in {1,2,3,4}{
\coordinate (A\i) at (-20+8*\i,6);
\begin{scope}[shift=(A\i)]

% positions

\coordinate (d1\i) at (0,0);
\coordinate (d2\i) at (0,3.5);
\coordinate (s1) at (-2,1);
\coordinate (s2) at (2,1);
\coordinate (x0) at (0,1);
\coordinate (x11) at (-1,1.5);
\coordinate (x12) at (-1,2.5);
\coordinate (x21) at (1,1.5);
\coordinate (x22) at (1,2.5);

%edges
\draw (d1\i)--(s1)--(x12)--(d2\i)--(x22)--(s2)--(d1\i);
\draw (d1\i)--(x0);
\draw (s2)--(x21)--(x0)--(x11)--(s1);
\draw
(x12)--(x0)--(x11)--(x22)--(x12)--(x21)--(d2\i)--(x11)--(x21)--(x0)--(x22);

%vertices
\draw (d1\i) node[blackvertex] {};
\draw (d2\i) node[blackvertex] {};
\draw (s1) node[vertex] {};
\draw (s2) node[vertex] {};
\draw (x0) node[vertex] {};
\draw (x11) node[vertex] {};
\draw (x12) node[vertex] {};
\draw (x21) node[vertex] {};
\draw (x22) node[vertex] {};

%labels
%\draw (d1\i) node[below right] {$a_1$};
%\draw (d2\i) node[above right] {$a_2$};
%\draw (s1) node[left] {$b_1$};
%\draw (s2) node[right] {$b_2$};
%\draw (x0) node[below right] {$z$};
%\draw (x11) node[below] {$x_1$};
%\draw (x12) node[above left] {$x_2$};
%\draw (x21) node[below] {$y_1$};
%\draw (x22) node[above right] {$y_2$};
\end{scope}
}

%%%% CLAUSES
% positions
\foreach \i in {0,1,...,\k}{
\coordinate (c\i) at (-90+\i*360/\l:3); }

\path (c0) ++(-0.5,-0.5) coordinate (x2) +(120:1) coordinate (x1);
\path (c0) ++(0.5,-0.5) coordinate (x3) +(60:1) coordinate (x4);

%edges
\draw (x1)--(x2)--(x3)--(x4);
\foreach \i in {1,2,...,\k}{
   \foreach \j in {1,2,...,\i}{
     \draw (c\i)--(c\j);
   }};
   \foreach \i in {1,2,...,\k}{
      \draw (x1)--(c\i)--(x4);
      }

%%%%%%%%%%%%%%%%%%%%%%%%%%%%%%%%%%%%%%%%%%%%%%%%%%%%%%%%%%%%%%%%%%%%%%%%%%%%%%%%%%%%%%
%                   FORMULA                       %
%%%%%%%%%%%%%%%%%%%%%%%%%%%%%%%%%%%%%%%%%%%%%%%%%%%%%%%%%%%%%%%%%%%%%%%%%%%%%%%%%%%%%%

\draw (c3)--(d11) (d21) .. controls +(+2,-0.5) .. (c3);
\draw (c2)--(d12) (d22) .. controls +(+2,-0.5) .. (c2);
\draw (c2)--(d13) (d23) .. controls +(-2,-0.5) .. (c2);
\draw (c1)--(d13) (d23) .. controls +(2,-0.5) and +(5.5,8) .. (c1);
\draw (c1)--(d14) (d24) .. controls +(-2,-0.5) .. (c1);

%%%%%%%%%%%%%%%%%%%%%%%%%%%%%%%%%%%%%%%%%%%%%%%%%%%%%%%%%%%%%%%%%%%%%%%%%%%%%%%%%%%%%%

%vertices
\foreach \i in {1,2,...,\k}{
\draw (c\i) node[clause] {};}
\draw (x1) node[vertex] {};
\draw (x2) node[vertex] {};
\draw (x3) node[vertex] {};
\draw (x4) node[vertex] {};

%labels
%\foreach \i in {1,2,...,\k}{
%\path (c\i) +(-30+\i*360/\l:0.4) node {$c_\i$};}

%\draw (x1) node[below left] {$p1$};
%\draw (x2) node[below] {$p2$};
%\draw (x3) node[below] {$p3$};
\draw (x4) node[below right] {$P$};

\draw (c1) node[below right] {$c_3$};
\draw (c2) node[right] {$c_2$};
\draw (c3) node[below left] {$c_1$};

\foreach \i in {1,2,3,4}{
\draw (d2\i) node[above] {$X_\i$};
}

\end{tikzpicture}

\end{center}
\caption{Example of the graph for formula $X_1 \wedge (X_2\vee
X_3)\wedge (X_3\vee X_4)$ }
\label{fig:vexample}
\end{figure}

Next we present a construction of the graph $G_{\cal F}$, when we are given a formula ${\cal F}$ with $k$ variables and $n$ clauses. We require that $k$ is even, otherwise we add a variable that appears in no clause.
For each variable $X$ in ${\cal F}$ we take a copy $W_X$ of the graph $W|\{a_1,a_2\}$ (that is, we assume that $a_1$ and $a_2$ in the copy $W_X$ are dominated in $G_{\cal F}$). For each disjunctive clause ${\cal C}_i$ in the formula we add a vertex $c_i$, and for each $X$ that appears in ${\cal C}_i$ we make $c_i$ adjacent to both $a_1$ and $a_2$ from the copy of $W_X$. Next, we add edges $c_ic_j$ between each two vertices corresponding to disjunctive clauses ${\cal C}_i,{\cal C}_j$ that appear in $\cal F$. Hence the vertices $c_i$, $1\le i\le n$, induce a clique $Q$ of size $n$. Finally, we add a copy $P:p_1p_2p_3p_4$ of a path $P_4$, and add edges $p_1c_i$ and $p_4c_i$ for $1\le i\le n$.
See Fig.~\ref{fig:vexample} for an example of the construction.

We call $W_X$ a  {\em widget subgraph} of $G_{\cal F}$, and in the notation for vertices in $W_X$ we add $X$
as an index to a vertex from $W_X$. For instance, a vertex that corresponds to $a_1$ in a widget subgraph $W_X$ will be denoted by $a_{1,X}$,
while the vertex that corresponds to $z$ in this subgraph will be denoted by $z_X$.
Denoting $A=\{a_{1,X},a_{2,X}\,|\,X \textrm{ variable in } {\cal F}\}$, we will consider the game played on $G_{\cal F}|A$.

%%%%%%%%%%%%%%%%%%%%%%%%%%%%%%%%%%%%%%%%%%%%%%%%%%%%%%%%%%%%%%%%%%%%%%%%%
% Gasper and Sandi believe that we do not need the following paragraph as we say it all later on whenever needed.
%
% In the proof of the main result, we often use the following moves and their interpretation, which connects the game played in a graph $G_{\cal F}$ with the POS-CNF game played in ${\cal F}$. Notably, if \D\ plays first in the widget subgraph $W_X$, then he plays $a_1$, by which each $c$, adjacent to it, is also dominated. This move of \D\ corresponds to the variable $X$ being set to TRUE by Player 1 in POS-CNF game. If \St\ plays first in the widget subgraph $W_X$ then she plays $b_1$, and this move of \St\ corresponds to the variable $X$ being set FALSE by Player~2 in POS-CNF game.
%%%%%%%%%%%%%%%%%%%%%%%%%%%%%%%%%%%%%%%%%%%%%%%%%%%%%%%%%%%%%%%%%%%%%%%%%

\medskip

The following observation will also be useful to consider the game played on $G_{\cal F}|A$.

\begin{figure}[ht!]
\begin{center}

\begin{tikzpicture}
[thick,scale=0.7,
     vertex/.style={circle,draw,inner sep=0pt,minimum
size=1.5mm,fill=white!100},
     blackvertex/.style={circle,draw,inner sep=0pt,minimum
size=1.5mm,fill=black!100},
     clause/.style={circle,draw,inner sep=0pt,minimum
size=1.5mm,fill=white!100}]

\def\k {4};
\def\l {5};

% positions
\foreach \i in {0,1,...,\k}{
\coordinate (c\i) at (-90+\i*360/\l:2); }

\path (c0) ++(-0.5,-0.7) coordinate (x2) +(120:1) coordinate (x1);
\path (c0) ++(0.5,-0.7) coordinate (x3) +(60:1) coordinate (x4);

%edges
\draw (x1)--(x2)--(x3)--(x4);
\foreach \i in {1,2,...,\k}{
   \foreach \j in {1,2,...,\i}{
     \draw (c\i)--(c\j);
   }}

\foreach \i in {1,2,...,\k}{
      \draw (x1)--(c\i)--(x4);
  %    % partial edges
   %   \draw (c\i) -- +(-60+360*\i/\l:0.5) (c\i)-- +(-120+360*\i/\l:0.5);
      }

%vertices
\foreach \i in {1,2,...,\k}{
\draw (c\i) node[clause] {};}

\draw (c1) node[blackvertex] {};
\draw (c2) node[blackvertex] {};
\draw (c4) node[blackvertex] {};
\draw (x1) node[vertex] {};
\draw (x2) node[vertex] {};
\draw (x3) node[vertex] {};
\draw (x4) node[vertex] {};

%labels
%\foreach \i in {1,2,...,\k}{
% \path (c\i) +(-30+\i*360/\l:0.4) node {$c_\i$};}
\draw (c1) node[right] {$c_1$};
\draw (c2) node[right] {$c_2$};
\draw (c3) node[left] {$c_3$};
\draw (c4) node[left] {$c_4$};
\draw (x1) node[left] {$p_1$};
\draw (x2) node[left] {$p_2$};
\draw (x3) node[right] {$p_3$};
\draw (x4) node[right] {$p_4$};
\end{tikzpicture}

\end{center}
\caption{The graph $H|\{c_1,c_2,c_4\}$ for $n=4$}
\label{fig:clause}
\end{figure}
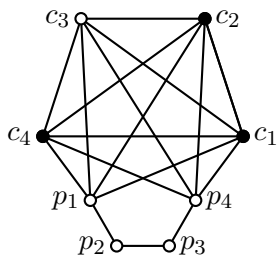

\begin{observation}
\label{obs:clauses}
Let $H$ be a graph isomorphic to the subgraph of $G_{\cal F}$ induced by the vertices from $Q\cup P$, and let $S\subseteq V(Q)$ be some vertices already dominated (see Fig.~\ref{fig:clause}).
\begin{enumerate}[(i)]
\item If $Q$ is not entirely dominated, that is, if $S\neq V(Q)$, then for any first move of \D\ in $H$, \St\ can reply and leave a vertex of $P$ undominated. Thus $\gamma_g(H|S)= 3$. 
\item If $S = V(Q)$, then $\gamma_g(H|S)=2$.
\item For any $S\subseteq V(Q)$, $\gamma_g'(H|S)=2$.
\end{enumerate}
\end{observation}

\begin{lemma}
\label{lem:D-win}
If Player 1 has a winning strategy for a formula ${\cal F}$ in the POS-CNF game, then \D\ can ensure at most $3k+2$ moves are played in the D-game played on $G_{\cal F}|A$ (i.e. $\gamma_g(G_{\cal F}|A) \le 3k+2$).
\end{lemma}

\proof
Assume that Player 1 has a winning strategy for ${\cal F}$ in the POS-CNF game, we give a strategy of \D\ that ensures at most $3k+2$ moves will be played in the D-Game on $G_{\cal F}|A$.
The general idea for \D\ is to ensure no more than three moves are played in each widget, while playing the vertex $a_{1,X}$ for each $X$ that would have been set TRUE in the POS-CNF game. This way he would dominate the vertices in $Q$. Then two moves suffice to dominate the vertices from $P$ and finish the game. 
 
To describe precisely the strategy of \D\ in $G_{\cal F}$, we use a simultaneously played POS-CNF game on $\cal F$. The first move of Dominator is to play $a_{1,X}$ where $X$ is the first variable Player~1 would set TRUE in the POS-CNF game.

\begin{itemize}
\item Whenever Staller makes a move, which is the first move in a widget subgraph $W_Y$, Dominator considers this move as Player~2 setting $Y$ FALSE in the POS-CNF game. Dominator then follows the POS-CNF winning strategy of Player~1; that is, he 
plays $a_{1,X}$ where $X$ is the next variable he would set TRUE in the POS-CNF game, as long as there are undefined variables. 

\item If Staller makes a move in a widget subgraph $W_X$ where one move was played already, Dominator answers in the same widget subgraph $W_X$ to guarantee that no more than three moves are played in $W_X$. 
\begin{itemize}
\item In particular, if $X$ is a variable set TRUE, then \D's second move in $W_X$ prevents the situation from Observation~\ref{obs:widget}(iii). 
\item If $X$ is a variable set FALSE, Dominator finishes to dominate $W_X$ ensuring that $a_{1,X}$ or $a_{2,X}$ is played, as noticed in Observation~\ref{obs:widget}(vi) (otherwise Staller could be allowed to play $a_{1,X}$ later on, if some adjacent vertices in $Q$ were still undominated\footnote{This is actually why such a complicated widget is needed for the variables, otherwise two adjacent vertices could serve as a widget subgraph.}).
\end{itemize}  
Note that possibly after \St's move, the widget subgraph is completely dominated already, which would prevent \D\ from playing in $W_X$. However, this would imply that only two moves are played in $W_X$. Then we claim \D\ can ensure at most $3k+2$ moves are played in total. He first plays in $Q$, ensuring at most three moves are played in $Q\cup P$ by Observation~\ref{obs:clauses}, then by Observation~\ref{obs:widget}(i) and (iv), he ensures at most three moves are played in every widget subgraph but $W_X$, that got only two moves. The total number of moves in the game is then at most $3k+2$.

\item If Staller played in $Q\cup P$, then using Observation~\ref{obs:clauses}(iii), Dominator finishes dominating $Q\cup P$ in the next move. By the previous paragraph, \D\ can ensure at most three moves are played in each widget subgraph, which implies that a total of at most $3k+2$ moves are played in the game, as desired. 
\end{itemize}
We may thus assume that Staller did not play in $Q \cup P$, and that all variables are set TRUE or FALSE. Since Dominator followed POS-CNF strategy, $\cal F$ is true and all vertices in $Q$ are thus dominated. Then Dominator can play $p_2$ and ensure no more than two moves are played in $Q\cup P$ (by Observation~\ref{obs:clauses}(ii)). Recalling that he also ensures at most three moves are played in each of the widget subgraphs (also in widget subgraphs that may not be entirely dominated yet), this implies $\gamma_g(G_{\cal F}|A) \le 3k+2$. 
\qed

\begin{lemma}
\label{lem:S-win}
If Player 2 has a winning strategy for a formula ${\cal F}$ in the POS-CNF game, then \St\ can ensure no less than $3k+3$ moves are played in the D-game played on $G_{\cal F}|A$ (i.e. $\gamma_g(G_{\cal F}|A) \ge 3k+3$).
\end{lemma}

\proof
Assume Player 2 has a strategy to win the POS-CNF game. We here describe a strategy for Staller to ensure at least $3k+3$ moves will be played in $G_{\cal F}|A$. The basic strategy of \St\ is to ensure that three moves will be played on each widget. This can be done with the following two rules:

\begin{itemize}
\item  Whenever \D\ makes a move in a widget subgraph $W_X$ in which no moves were made before, \St\ responds with $b_{1,X}$ or $b_{2,X}$ in the same widget subgraph $W_X$, ensuring at least three moves are played in $W_X$ (this can be done due to Observation~\ref{obs:widget}(i)). 
\item Whenever \St\ makes a move in a widget where no moves were played before, she plays $b_{1,X}$ thus ensuring three moves will be played (see~Observation~\ref{obs:widget}(iv)).
\end{itemize}

Note that two additional moves will be required to dominate $P$, so if \St\ can ensure either a third move in $Q\cup P$ or a fourth move in any widget subgraph $W_X$, we get $3k+3$ moves in the whole game, as desired.

%%In addition, note that to dominate the path $P$, at least two moves are needed during %the game, hence \St\ easily ensures at least $3k+2$ moves in the whole game. If %Staller manages to force an extra move in any of the mentioned subgraphs, we get our %result.
%Suppose Dominator plays a move in some widget subgraph $W_X$, which is the third move played in $W_X$, and which makes $W_X$ dominated. If there exists another widget subgraph $W_Y$ in which two moves were played and two vertices of $W_Y$ remain undominated (in particular, this happens when the two moves are played by vertices $a_{1,Y}$ and $b_{1,Y}$) then, by Observation~\ref{obs:widget}(iii), Staller can enforce four moves in that widget subgraph. Since in all other widget subgraphs at least three moves will be played, and at least two moves must be played in $Q\cup P$, we infer that at least $3k+3$ moves will be played in the game. Since the proof is done in this case, we may assume that \D\ does not play in this way. We resume this as follows.
%
%\begin{remark} 
%\label{rem2}
%We may assume that if \D\ played a third move in some widget subgraph $W_X$, by which $W_X$ becomes dominated, there is no other widget subgraph $W_Y$ in which exactly two moves were played and two vertices of $W_Y$ remain undominated. 
%\end{remark}

To describe the rest of \St's strategy in $G_{\cal F}$, we use a simultaneously played POS-CNF game on $\cal F$. Whenever \D\ plays $a_{1,X}$ or $a_{2,X}$ in a widget subgraph $W_X$, we transpose it to Player~1 setting the variable $X$ TRUE in the POS-CNF game, and when \St\  plays the third move without $a_{1,X}$ nor $a_{2,X}$ being played, or when she plays first in a widget subgraph $W_X$, this is transposed to Player 2 setting $X$ FALSE. Other \D's moves and their responses should not impact the POS-CNF game and are called \emph{neutral sequences}.

We now describe \St's response to the various moves of \D. \St's strategy mimics Player~2 strategy in the POS-CNF game and therefore, the clause vertices of $Q$ are not all dominated until some vertex is played in $Q\cup P$.

\begin{itemize}
 
\item Whenever \D\ plays the first move of some widget subgraph $W_X$, she responds (according to her basic strategy) by playing $b_{1,X}$ or $b_{2,X}$ on $W_X$ in such a way that not all vertices of $W_X$ are dominated. In the case \D\ played $a_{1,X}$ or $a_{2,X}$, \St\ considers Player~1 set the variable $X$ TRUE in the POS-CNF game on $\cal F$. (Note that it is \D's turn again and he could make a move that would set another variable TRUE in the corresponding POS-CNF game, while Player~2 did not set a variable FALSE. Then we do not try to mimic the POS-CNF game anymore; we deal with this case at the end of the proof.) If \D\ plays another vertex as the first move in $W_X$, this is a typical neutral sequence. 

\item \D\ plays the third move in a widget subgraph $W_X$, by which $W_X$ becomes dominated. If this move is $a_{1,X}$ or $a_{2,X}$ and none of these was played before, this corresponds to setting the variable TRUE. Alternatively, $a_{1,X}$ or $a_{2,X}$ have been played by \D\ as the first move in $W_X$, and by the third move in $W_X$ \D\ prevents four moves will be played in $W_X$, see Observation~\ref{obs:widget}(iii). In this case it is Player~2's turn to play in the POS-CNF game, and \St\ responds by selecting a variable $Y$ according to the strategy in the POS-CNF game (there is always an unset variable since the number of variables to be set is even). Either no vertices were yet played in $W_Y$, then \St\ set $Y$ FALSE by playing $b_{1,Y}$ as the first move in $W_Y$, or \D\ played a neutral move in $W_Y$ and \St\ answered already, then \St\ sets the variable $Y$ FALSE by playing $b_{1,Y}$, $b_{2,Y}$ or $z_Y$, whichever finishes dominating the widget $W_Y$.

\item \D\ plays in some widget subgraph $W_Y$, where $Y$ has been set FALSE by Staller (playing $b_{1,Y}$). If this is the fourth move in $W_Y$, by our first argument, we get that at least $3k+3$ moves are played in the game as desired. Suppose thus it is the second move. By Observation~\ref{obs:widget}(v), \D\ playing $a_{1,Y}$ or $a_{2,Y}$ implies that \St\ can enforce four moves on $W_Y$, which then implies the total of $3k+3$ moves and conclude the proof. For any other move of \D\ (e.g. $z_{Y}$ as suggested by Observation~\ref{obs:widget}(iv)), \St\ replies $b_{2,Y}$ or $z_Y$, making $W_Y$ dominated with the third move. This is also a neutral sequence.

\item \D\ plays in $Q\cup P$. By the above, the moves of \D\ and \St\ were emulating the moves of the POS-CNF game, in which Player~2 has a winning strategy. This implies that not all vertices of $Q$ are dominated at this point, even if all variables were set. Thus, by Observation~\ref{obs:clauses}(i), Staller can ensure three moves are played in $Q\cup P$ and this makes a total of $3k+3$ moves.
\end{itemize}

We now consider the case we postponed, when \D\ tries to make a sequence of at least two moves that correspond to setting a variable TRUE in the POS-CNF game while \St\ did not have the opportunity to set a variable FALSE. We thus are in the situation when \D\ played  $a_{1,X}$ or $a_{2,X}$ as the first move in a widget subgraph $W_X$, forcing an answer of \St\ on $b_{1,X}$. Then, possibly after some neutral sequences, \D\ sets another variable TRUE. If he does so by playing a third move in a widget subgraph $W_Y$ where a neutral sequence was played before, then, following Observation~\ref{obs:widget}(iii), \St\ replies in $W_X$ to enforce four moves in that widget and we get at least $3k+3$ moves in total. Otherwise \D\ set another variable TRUE playing $a_{1,Y}$ or $a_{2,Y}$ as the first move in another widget subgraph $W_Y$. Then \St\ must reply in $W_Y$ to ensure three moves will be played in that widget, but we claim she can enforce four moves either in $W_X$ or in $W_Y$. Indeed, whenever \D\ will make a move which is not the first move in a widget subgraph (it may be a third move that finishes dominating $W_X$ or $W_Y$ but not both), then \St\ can ensure four moves are played in $W_Y$ or $W_X$. This concludes the proof.
\qed

From Lemma~\ref{lem:D-win} and \ref{lem:S-win}, we can infer the following.

\begin{corollary}
\label{cor:main}
Player 1 has a winning strategy for a formula ${\cal F}$ in the POS-CNF game if and only if $\gamma_g(G_{\cal F}|A)\le 3k+2$.
\end{corollary}

We now use the reduction from the POS-CNF game to {\sc Game Domination Problem} as presented in Lemmas~\ref{lem:D-win} and~\ref{lem:S-win} to prove our main theorem. In particular, we show that having a partially dominated graph is not necessary.

\begin{theorem}
{\sc Game Domination Problem} is log-complete in PSPACE.
\label{t:main}
\end{theorem}
\proof
 Let us describe how to build a graph where no vertices are assumed already dominated. 

Consider a formula $\cal F$. First add a variable $X_0$ and modify $\cal F$ into a formula $\cal F'$ in such a way that variable $X_0$ is inserted into every clause.
Now take the corresponding graph $G_{\cal F'}|A$ where $A$ contains the vertices $a_{1,X}, a_{2,X}$ for all variables $X$ (including $X_0$). Add to $G_{\cal F'}$ (where $A$ is no longer dominated) a star $K_{1,3}$ with center $v$, and add an edge between $v$ and each vertex in $A$; we denote the resulting graph by $G'_{\cal F}$. 
The number of all variables (together with $X_0$) is $k+1$. 

{\bf Claim.} $\gamma_g(G'_{\cal F})\le 3k+6$ if and only if Player~1 wins the POS-CNF game on $\cal F$. 

{\bf Proof of Claim.} First we prove that $v$ is an optimal first move of \D, and
an optimal response of \St\ is $b_{1,X_0}$. 

Suppose that \D\ does not play $v$ as his first move. Then the response of \St\ is to play one of the leaves adjacent to $v$. If the next move of \D\ is to play $v$, then by the basic strategy of \St\ as presented in the proof of Lemma~\ref{lem:S-win} there will still be at least three moves played in each of the widget subgraphs and at least two moves in $Q\cup P$. Altogether this yields $2+3(k+1)+2=3k+7$ moves in the game.
If, on the other hand, the second move of \D\ is again not playing $v$, \St\ again responds by playing an undominated leaf. In this way, \D\ could gain an advantage by dominating one of the widget subgraphs in only two moves, however the advantage is lost by \St's second move on the leaf, which again yields $3k+7$ moves, if \D's third move is playing $v$. Otherwise, if \D's third move is not playing $v$, \St\ goes with her basic strategy. She can ensure at most one widget subgraph ends up with only two moves played, and a third move is required to dominate the third leaf. Then the total number of moves is $3k+2$ on the widget subgraphs, two for $P$ and three on the $K_{1,3}$, which makes $3k+7$. We derive that if \D\ does not play $v$ as his first move, there is no way that the number of moves would be less than $3k+7$. 

So we may assume that \D\ plays $v$ as his first move. If Staller does not respond by playing $b_{1,X_0}$, then in the next move \D\ can play $a_{1,X_0}$, by which $Q$ is dominated. Using the same strategy as in Lemma~\ref{lem:D-win}, \D\ can ensure that at most three moves are played in each widget subgraph, and since $Q$ is dominated, two moves suffice for $Q\cup P$. Altogether this sums up to at most $3k+6$ moves, which means that \St's second move must be $b_{1,X_0}$, if she wishes to have a chance to enforce more than $3k+6$ moves. 

Now, from $d_1=v,s_1=b_{1,X_0}$, we derive that all vertices of $A$ are dominated, that the (artificial) variable $X_0$ is set FALSE, and two additional moves in $W_{X_0}$ will be played, none of which is $a_{1,X_0}$ nor $a_{2,X_0}$ (if \D\ played one of these two vertices, \St\ can enforce four moves on the widget subgraph, by Observation~\ref{obs:widget}(v)). The rest of the game thus corresponds to the game as in the proof of Lemmas~\ref{lem:D-win} and~\ref{lem:S-win}, in which there are at most $3k+2$ moves if and only if Player~1 wins the POS-CNF game. Together with four additional moves (one playing $v$ and three in $W_{X_0}$), we deduce that $\gamma_g(G'_{\cal F})\le 3k+6$ if and only if Player~1 wins the POS-CNF game on $\cal F$. 

{\bf \qed (Claim)}

Finally, observe that the reduction from POS-CNF to {\sc Game Domination Problem} can be computed with a working space of size $O(\log(k+n))$; giving an explicit algorithm is a routine work. Recalling that POS-CNF is log-complete in PSPACE~\cite{ts-78} completes the proof.
\qed

%%%%%%%%%%%%%%%%%%%%%%%%%%%%%%%%%%%%%%%%%%%%%%%%%%%
\section{Conclusions and open problems}
\label{sec:conclude}
%%%%%%%%%%%%%%%%%%%%%%%%%%%%%%%%%%%%%%%%%%%%%%%%%%%
We can modify the reduction from the proof to apply it for the PSPACE-completeness of the Staller-start game domination problem, which can be done in a similar way as in the proof of Theorem~\ref{t:main}. Given a formula $\cal F$, we consider the formula $\cal F'$, where we add a variable $X_0$ that we insert into every clause. Clearly, if Player~2 starts the POS-CNF game on $\cal F'$, he must set variable $X_0$ FALSE as his first move to have a chance to win. Then, the winner of the POS-CNF game on the formula $\cal F'$ where Player~2 starts, is the same as the winner of the POS-CNF game on $\cal F$. Using this knowledge, it is straightforward to use the above reduction for the Staller-start domination game.

A related problem is to try to find families of graphs where the game domination number can be computed efficiently. Such families obviously include those where the exact formula is known (as already mentioned, these include paths, cycles and combs). Another family where we expect the domination game number can be computed in polynomial time is the class of proper interval graphs. For these graphs, it looks like both players' strategy can be described by a greedy algorithm, though we did not manage to prove it. Hence we pose:

\begin{question} \label{prob:interval}
Can the game domination number of (proper) interval graphs be computed in polynomial time?
\end{question}

In particular, interval graphs are also dually chordal graphs, which are proven in~\cite{doko-2015} to be the so-called no-minus graphs (that is, the graphs $G$ for which for all subsets of vertices $U$, $\gamma_g(G|U)\le\gamma_g'(G|U)$). In that paper,  stronger relations between the game domination number of the disjoint union of two graphs and the game domination number of the components are given for no-minus graphs, which could prove useful in Question~\ref{prob:interval}. On the other hand, it seems likely that the decision problem remains PSPACE-complete even when restricted to split graphs, which are also proven to be no-minus in~\cite{doko-2015}. Note that this would be a similar dichotomy as proved in~\cite{heggernes-2012} for the {\sc Role Assignment} problem which can be solved in polynomial time on proper interval graph and is {\sc Graph Isomorphism}-hard on chordal graphs.

Observe that the domination game on split graphs transposes to a game on hypergraphs where Dominator chooses edges and Staller chooses vertices (not chosen before nor belonging to a chosen edge), and the game ends after all vertices either are chosen or belong to a chosen edge. The aim of Dominator is again to finish the game as soon as possible, while Staller tries to have the game last for as long as possible. The game where both players must choose hyperedges may be of independent interest.

\section*{Aknowledgements}

% The authors wish to thank anonymous referees for their comments that helped to improve the presentation of the paper. 
This work was done in the frame of the bilateral France-Slovenian project BI-FR/13-14-PROTEUS-003 entitled Graph domination, supported in part by the French Minist\`ere des Affaires \'etrang\`eres (MAE) and Minist\`ere de l’Enseignement sup\' erieur et de la Recherche (MENESR) and by the Slovenian Research Agency (ARRS). B.B. and S.K. are supported in part by the Ministry of Science of Slovenia under the grant P1-0297. B.B., S.K. and G.K. were also supported in part by the Ministry of Science of Slovenia under
the grants P1-0297 and N1-0043. P.D. and G.R. got financial support from the French State, managed by the French National Research Agency (ANR) in the frame of the ``Investments for the future'' Programme IdEx Bordeaux - CPU (ANR-10-IDEX-03-02).

%%%%%%%%%%%%%%%%%%%%%%%%%%%%%%%%%%%%%%%%%%%%%%%%%%%%%%%%%%%%%%%%%

\end{document}